\documentclass[12pt]{amsart}
\usepackage{amsmath,amstext,amsfonts,amsthm,amssymb}
\usepackage{eucal}
\usepackage{graphicx}
\renewcommand{\subsubsection}[1]{\addtocounter{subsubsection}{1}
{\ \\[3pt]\bf \thesubsubsection. \  #1} }

\swapnumbers

{  \theoremstyle{definition}

}
\newcommand{\iso}{\overset{\sim}{\longrightarrow}}

\newcommand{\lra}{\longrightarrow}

\newcommand{\bea}{\begin{eqnarray*}}
\newcommand{\eea}{\end{eqnarray*}}
\newcommand{\bean}{\begin{eqnarray}}
\newcommand{\eean}{\end{eqnarray}}

\newcommand{\blambda}{{\bar{\lambda}}}

\newcommand{\CG}{\mathcal{G}}

\newcommand{\BC}{\mathbb{C}}

\newcommand{\nc}{\newcommand}

\nc{\Id}{\text{Id}}
\nc{\la}{\lambda}

\begin{document}


\centerline{\bf INVARIANT FUNCTIONALS}

\bigskip\bigskip

\centerline{\bf AND ZAMOLODCHIKOVS' INTEGRAL}

\bigskip\bigskip


\centerline{Bui Van Binh and Vadim Schechtman}

\bigskip\bigskip




\bigskip\bigskip

Consider a triple complex integral
$$
I(\sigma_1,\sigma_2,\sigma_3) = 
\int_{\BC^3} (1 + |z_1|^2)^{-2\sigma_1}(1 + |z_2|^2)^{-2\sigma_2}(1 + |z_3|^2)^{-2\sigma_3}
$$
$$
|z_1 - z_2|^{2\nu_3-2}|z_2 - z_3|^{2\nu_1-2}|z_3 - z_1|^{2\nu_2-2} dz_1dz_2dz_3.
\eqno{(1)}
$$
It 
has appeared in a remarkable paper [ZZ] in connection with the Liouville model of the conformal field theory. Here $\sigma_i\in \BC$, 
$$
\nu_i = \sigma_{i+1} + \sigma_{i+2} - \sigma_i,\ i\mod\ 3,
\eqno{(2)}
$$ 
$dz := dx dy,\ z = x + iy$, stands for the standard Lebesgue measure on $\BC$. 
The integral converges if $\Re \sigma_i > 1/2,\ \Re \nu_i > 0$, $i = 1, 2, 3$. 
The aim of this note is to prove the following assertion. 

{\bf Theorem 1.}
$$
I(\sigma_1,\sigma_2,\sigma_3) = \pi^3\frac{\Gamma(\nu_1 + \nu_2 + \nu_3 - 1)
\Gamma(\nu_1)\Gamma(\nu_2)\Gamma(\nu_3)}{\Gamma(2\sigma_1)\Gamma(2\sigma_2)\Gamma(2\sigma_3)}
\eqno{(3)}
$$

This formula is given without proof in [ZZ]. A proof of (3) has appeared in [HMW]; it uses some complicated  
change of variables. A different proof using 
Macdonald type constant term identities may be found in [BS]. In this note we provide still 
another proof of (3) which uses an elegant method of Bernstein and Reznikov who computed a similar real integral, cf. [BR]. Their procedure 
is based on some multiplicity one considerations which allow one to 
reduce the computation of (1) to a computation of a Gaussian integral, cf. Thm 2 below; 
this last integral is of independent interest. 

{\it Notation.} For a smooth variety $Y$ we denote by $C(Y)$ the space of $C^\infty$ functions 
$f: Y\lra \BC$.

{\it Principal series.} Let $G := SL_2(\BC)\supset K := SU(2)$; $X = \BC^2\setminus \{0\}$. Let 
$C(X)$ denote the space of smooth functions $f:\ X\lra \BC$; for $\lambda\in \BC$ let 
$C_\lambda(X)\subset C(X)$ be the subspace of $f$ such that $f(ax) = |a|^{2\lambda}f(x)$ for all 
$a\in\BC^*$. The group $G$ is acting on $C_\lambda(X)$ by the rule $(gf)(x) = f(xg),\ g\in G$; 
we denote this representation by $V_\lambda$. 

{\it Line realization}, cf. [GGPS]. The map $f\mapsto f(x,1)$ induces an isomorphism 
$\phi: V_\lambda \iso V'_\lambda$ where $V'_\lambda$ is a subspace of $C(\BC)$ depending on $\lambda$. 

The induced action of $G$ on $V'_\lambda$ is given by 
\newline $(gf)(x) = |bx + d|^{2\lambda}f((ax+c)/(bx+d)),\ 
g = \left(\begin{matrix}a & b\\c & d\end{matrix}\right)$.  

{\bf Lemma 1.} {\it 
The dimension of the space $Hom_K(V_\lambda,\BC)$ of continuous $K$-invariant maps 
$V_\lambda\lra\BC$ is equal to $1$.}

{\bf Proof.} Let $\ell\in Hom_K(V_\lambda,\BC)$. Choose a Haar measure on $K$ such that 
$\int_K dk = 1$. Since $\ell$ is continuous and $K$ invariant, for all $f\in V_\lambda$ 
$\ell(f) = \ell(\bar f)$ where $\bar f(x) = \int_K f(kx) dk$. But $\bar f$ is homogeneous 
and $K$-invariant, and such a function is unique up to a multiplicative constant. 
$\square$

{\bf Lemma 2.} (i) {\it The functional $\ell'_\lambda:\ V'_{\lambda-2}\lra\BC$ given by 
$$
\ell'_\lambda(f) = \int_\BC (1 + |z|^2)^{-\lambda} f(z) dz
$$
 is $K$-invariant.The integral converges for all $\lambda\in \BC,\ f\in V'_{\lambda-2}$. 
$\ell'_\lambda$ is $G$-invariant iff $\lambda = 0$.}    

(ii) {\it The functional $\ell_\lambda:\ V_{\lambda-2}\lra\BC$ given by 
$$
\ell_\lambda(f) = \frac{1}{\pi}\int_{\BC^2} f(z_1,z_2) e^{-|z_1|^2 - |z_2|^2} dz_1dz_2
$$ 
is $K$-invariant. The integral converges for $\Re \lambda > 0$.} $\square$

By Lemma 1, $\ell_\lambda = c_\lambda \ell'_\lambda\circ \phi$ for some $c_\lambda\in\BC$.

{\bf Lemma 3.} $c_\lambda = \Gamma(\lambda)\pi^{-1}$. $\square$

More generally, given $n$ complex numbers $\blambda =(\lambda_1, \ldots, \lambda_n)$ we define 
a $G$-module 
$$
V_{\blambda} = \{f\in C(X^n)| f(a_1x_1,\ldots,a_nx_n) = \prod_{i=1}^n\ |a_i|^{2\lambda_i}f(x_1,\ldots,x_n)\} 
$$ 
with the diagonal action of $G$. 

The map $f\mapsto f((x_1,1),\ldots,f(x_n,1))$ induces an isomorphism of $G$-modules 
$$
\phi:\ V_{\blambda}\iso V'_{\blambda} \subset C(\BC^n)
$$  
with the action of $G$ on $V'_{\blambda}$ given by 
$$
(gf)(x_1,\ldots,x_n) = \prod_{i=1}^n |bx_i + d|^{-2\lambda_i}
f\biggl(\frac{ax_1 + c}{bx_1 + d},\ldots,\frac{ax_n + c}{bx_n + d}\biggr),\ 
g = \left(\begin{matrix}a & b\\c & d\end{matrix}\right).
$$

We have $K$-invariant functionals $\ell'_{\blambda}:\ V'_{\lambda_1 - 2,\ldots,\lambda_n-2}\lra \BC$, 
$$
\ell'_{\blambda}(f) = \int_{\BC^n} \prod_{i=1}^n (1 + |z_i|^2)^{-\lambda_i}f(z_1,\ldots,z_n) dz_1\ldots dz_n, 
$$
(the integral converges for all $\lambda_i$)

and $\ell_{\blambda}:\ V_{\lambda_1 - 2,\ldots,\lambda_n-2}\lra \BC$,
$$
\ell_{\blambda}(f) = \int_{(\BC^2)^n} f(x_1,y_1,\ldots,x_n,y_n)e^{- \sum_{i=1}^n (|x_i|^2 + |y_i|^2)}dx_1dy_1\ldots dx_ndy_n
$$
(the integral converges if $\Re\lambda_i > 0$ for all $i$). 

{\bf Corollary 1.} 
$$
\ell_{\blambda} = \pi^{-n}\prod_{i=1}^n \Gamma(\lambda_i) \cdot \ell'_{\blambda}\circ\phi. 
$$

Indeed, the factorisable functions 
$\prod f_i(x_i),\ f_i\in V_{\lambda_i}$ are dense in $V_{\lambda_1, \ldots, \lambda_n}$, 
and we apply the Fubini theorem. 

For a function $f:\ Y\lra \BC$ where $Y\subset \BC^n$ with the complement of measure zero we denote 
$$
\CG(f) = \frac{1}{\pi^{n/2}}\int_{\BC^n} f(z_1,\ldots,z_n)e^{-\sum_{i=1}^n |z_i|^2} dz_1\ldots dz_n
$$
(when the integral converges). The following proposition 
is a complex analog of properties of real Gaussian integrals from [BR]. 

{\bf Proposition 1.} (i) {\it Set $r(z) = (\sum_{i=1}^n |z_i|^2)^{1/2}$. Then $\CG(r^s) = 
\Gamma(s/2 + n)/\Gamma(n)$.}

(ii) {\it If $h(z) = \sum c_iz_i$ then $\CG(|h|^s) = \Gamma(s/2+1)||f||^s$ where 
$||f|| = (\sum_{i=1}^n |c_i|^2)^{1/2}$.}

(iii) {\it Consider the space $\BC^4 = (\BC^2)^2$ and the determinant 
$$
d(w_1,w_2) =  z_{11}z_{22} - z_{12}z_{21},\ w_i = (z_{1i}, z_{2i}). 
$$  
Then $\CG(|d|^s) = 
\Gamma(s/2+1)\Gamma(s/2+2)$.} $\square$ 

Now consider the space $L = (\BC^2)^3$ whose elements we shall denote 
$(w_1,w_2,w_3),\ w_i \in \BC^2$. Consider the following function on $L$:
$$
K_{\nu_1,\nu_2,\nu_3}(w_1,w_2,w_3) = |d(w_1,w_2)|^{2\nu_3-2}|d(w_1,w_3)|^{2\nu_2-2}|d(w_2,w_3)|^{2\nu_1-2}
$$

{\bf Theorem 2.} 
$$
\CG(K_{\nu_1,\nu_2,\nu_3}) = \Gamma(\nu_1 + \nu_2 + \nu_3 - 1)
\Gamma(\nu_1)\Gamma(\nu_2)\Gamma(\nu_3)
$$

The proof goes along the same lines as in [BR],  where a similar real 
Gaussian integral is calculated. 

Now we easily deduce Theorem 1. The function $K_{\nu_1,\nu_2,\nu_3}$ belongs 
to the space $V_{2\sigma_1 - 2, 2\sigma_2 - 2, 2\sigma_3 - 2}$ where 
$\sigma_i$ are defined from (2). We have 
$$
\CG(K_{\nu_1,\nu_2,\nu_3}) = \ell_{2\sigma_1, 2\sigma_2, 2\sigma_3}(K_{\nu_1,\nu_2,\nu_3}).
$$ 
On the other hand, by Cor. 1 
$$
I(\sigma_1, \sigma_2, \sigma_3) = \ell'_{2\sigma_1, 2\sigma_2, 2\sigma_3}(K'_{\nu_1,\nu_2,\nu_3}) = 
\pi^3\prod_{i=1}^3 \Gamma(2\sigma_i)^{-1}\ell_{2\sigma_1, 2\sigma_2, 2\sigma_3}(K_{\nu_1,\nu_2,\nu_3})
$$
where $K'_{\nu_1,\nu_2,\nu_3} = \phi(K_{\nu_1,\nu_2,\nu_3})$, 
This, together with Thm 2, implies Thm 1. $\square$

The functional $\ell: V_{-2\sigma_1,-2\sigma_2,-2\sigma_3} \lra \BC$, 
$\ell(f) = \int_{\BC^3} f K'_{\nu_1,\nu_2,\nu_3}dz_1dz_2dz_3$, 
is $G$-invariant. Similarly to the real case treated in [BR], the space of such $G$-invariant functionals is one-dimensional, cf. 
[L]. This is a "methaphysical reason", why the integral (1) (the value of this functional on some spherical vector) is computable in terms of 
$\Gamma$ functions.  

We are grateful to A.Reznikov for discussions and for the notes by A.Yomdin about a $p$-adic 
analogue of the Bernstein - Reznikov method which were very useful for us.

 \bigskip\bigskip

\centerline{\bf References}

\bigskip\bigskip

[BR] J.Bernstein, A.Reznikov, Estimates of automorphic functions, {\it Moscow Math. J.} 
{\bf 4} (2004), 19 - 37. 

[BS] Bui V.B., V.Schechtman, Remarks on a triple integral, {\it Moscow Math. J.}, to appear; arXiv:1204.2117.

[GGPS] I.M.Gelfand, I.M.Graev, I.I.Pyatetsky-Shapiro, Generalized functions, Vol. {\bf 5}.  

[HMW] D.Harlow, J.Maltz, E.Witten, Analytic continuation of Liouville theory, 
arXiv:1108.441.

[L] H.Loke, Trilinear forms of $\frak{gl}_2$, {\it Pacific J. Math.} {\bf 197} (2001), 119 - 144. 

[ZZ] A.Zamolodchikov, Al.Zamolodchikov, Conformal bootstrap in Liouville 
field theory, {\it Nucl. Phys.} {\bf B477} (1996), 577 - 605.

\bigskip\bigskip

Institut de Math\'ematiques de Toulouse, 
Universit\'e P.Sabatier, 31062 Toulouse, France

\end{document}